\documentclass[psamsfont]{amsart}

\usepackage{graphicx}
\usepackage{amssymb}
\usepackage{xypic}
\usepackage{hyperref}
\usepackage{psfrag}
\usepackage{amsmath,amssymb}
\usepackage{dsfont}\let\mathbb\mathds
\usepackage{pstricks}
\usepackage{subfigure}

\numberwithin{equation}{section}
\newtheorem{theorem}[equation]{Theorem}
\newtheorem{lemma}[equation]{Lemma}
\newtheorem{corollary}[equation]{Corollary}

\newtheorem{conjecture}{Conjecture}

\theoremstyle{remark}
\newtheorem{rem}[equation]{Remark}

\theoremstyle{definition}
\newtheorem{definition}[equation]{Definition}

\renewcommand{\qedsymbol}{\nobreak \ifvmode \relax \else \ifdim\lastskip<1.5em
  \hskip-\lastskip \hskip1.5em plus0em minus0.5em \fi \nobreak \vrule
  height0.75em width0.5em depth0.25em\fi}

\DeclareMathOperator{\grad}{\nabla}

\DeclareMathOperator{\lap}{\Delta}

\newcommand{\abs}[1]{\lvert #1 \rvert}

\newcommand{\nts}[1]{\marginpar{#1}}
\renewcommand{\nts}[1]{}
\newcommand{\R}{\mathbf{R}}

\newcommand{\N}{\mathbf{N}}

\newcommand{\ball}{\mathcal{B}}

\def\del{\partial}              

\def\bR{\mathbb R}          

\def\mF{\mathcal{F}}            
\def\mN{\mathcal{N}}

\begin{document}

\title[Extrema of low eigenvalues of the D--N Laplacian on a disk]{Extrema of low eigenvalues of the Dirichlet--Neumann Laplacian on a disk}
\author{Eveline Legendre}
\address{E.L.: Universit\'e de Montr\'eal\\
CP 6128 Succ. Centre ville\\
Montr\'eal, QC H3C 3J7\\
Canada}
\email{legendre@dms.umontreal.ca}
\thanks{This paper is based on the author's M.Sc. thesis, being carried out under the supervision of Professor Iosif Polterovich at Universit\'e de Montr\'eal.}

\maketitle

\begin{abstract}
We study extrema of the first and the second mixed eigenvalues of the Laplacian on the disk among some families of Dirichlet--Neumann boundary conditions. We show that the minimizer of the second eigenvalue among all mixed boundary conditions lies in a compact $1$--parameter family for which an explicit description is given. Moreover, we prove that among all partitions of the boundary with bounded number of parts on which Dirichlet and Neumann conditions are imposed alternately, the first eigenvalue is maximized by the uniformly distributed partition.
\end{abstract}

\section{Introduction}
\label{intro}

Let $\Omega$ be a $d$--dimensional Euclidean domain with Lipschitz boundary $\partial \Omega$. Let $\partial_D
\Omega$ be an open subset of $\partial \Omega$, and let $\partial_N \Omega$ be the
open remainder $\partial \Omega \setminus \overline{\partial_D \Omega}$. Consider the mixed Dirichlet--Neumann eigenvalue problem determined by imposing
Dirichlet condition on $\partial_D \Omega$ and Neumann condition on
$\partial_N \Omega$ for the eigenfunction of the Laplacian. That is
\begin{equation}
  \begin{split}
    -\lap u &= \lambda(\Omega, \partial_D \Omega) u,\\
    \left. u \right|_{\partial_D \Omega} &= 0,\\
    \left. \frac{\partial u}{\partial n} \right|_{\partial_N \Omega} &= 0\\
  \end{split}
\end{equation}
where $\frac{\partial u}{\partial n}$ denotes the derivative of $u$ with respect to the vector field, normal to the boundary of $\Omega$. So there is a infinite, discrete and positive spectrum which we ordinate and denote $\lambda_1(\Omega, \del_D \Omega) < \lambda_2(\Omega, \del_D \Omega) \leq \lambda_3(\Omega, \del_D \Omega) \dots $   This boundary value problem is also called Zaremba's problem, see~\cite{z:zaremba}.

Our interest is in understanding the dependence of these eigenvalues on the geometric properties of the partition of the boundary into Dirichlet and Neumann parts. In order to do so we study the extrema, minimum or maximum, of low eigenvalues among all parts of the boundary $\partial_D \Omega$ with fixed volume. To the author knowledge this point of view was first presented in Denzler's works on the first mixed eigenvalue, see~\cite{denzler:bounds},~\cite{denzler:disc}.

Mixed Dirichlet--Neumann problems appear naturally in many physical models and the low eigenvalues ($\lambda_1(\Omega, \del_D \Omega), \lambda_2(\Omega, \del_D \Omega)\; \dots $) are especially important in these models since they correspond to lower energy. The following interpretation, first proposed by Walter Craig as reported by Denzler \cite{denzler:bounds}, can be particularly enlightening: Interpreting the domain $\Omega$ as a room, Neumann's parts $\del_N \Omega$ may be viewed as perfectly insulated walls while Dirichlet's parts $\del_D \Omega$ are non-insulated windows. From this point of view, low eigenvalues (depending on initial distribution of temperature) determine generically the rate of heat diffusion through windows after large time (here we ignore convection). Hence, loosely speaking, extrema of low eigenvalues can be understood as extrema of heat loss through windows after large time.

Extrema of the first mixed eigenvalue have been studied in different cases. But extrema of higher mixed eigenvalues are essentially unknown at this time. This paper focuses on the study of the minimal arrangement of boundary conditions for the second eigenvalue on the disk. More precisely, we point out an explicit, compact, $1$--parameter family of boundary conditions containing the minimizer of the second mixed eigenvalue.

Note that it is natural to minimize, or maximize, eigenvalues over the family of Dirichlet's parts of given length since $\partial_{D_1} \ball \subset \partial_{D_2} \ball$ implies $\lambda_k(\partial_{D_1} \ball) \leq \lambda_k(\partial_{D_2} \ball)$. We define:

\begin{definition}[Boundary family $\mF_\ell$]
  For $ 0 < \ell < 2\pi$, the space of Dirichlet's
  part, $\partial_D \ball \subset \partial \ball$, with length $\ell$ as,
  \begin{multline}
    \mF_{\ell} := \{\partial_D \ball \subset \partial\ball \; : \; |\partial_D
    \ball|=\ell\ \\\text{ and } \partial_D \ball \text{ has a finite number of
      connected components}\}.
  \end{multline}
\end{definition}

In addition, we define the following subfamily of special arrangements of the boundary conditions which is, as we shall see, of particular interest.

\begin{definition}[Uniform $n$--partition]
  \label{unifKpart}
  For $n \in \N$, a uniform $n$--partition of length $\ell$, denoted by
  $\Gamma_n\in\mF_\ell$, is the union of $n$ connected parts of equal length
  $\frac{\ell}{n}$ uniformly distributed in $\partial \ball$.
\end{definition}

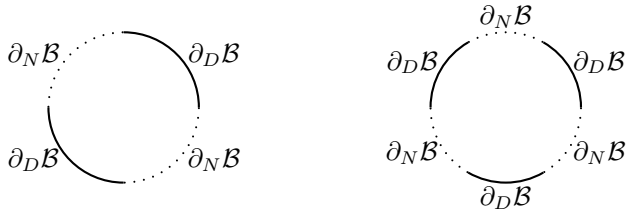
\begin{figure}[h]
\centering
\begin{minipage}{10cm}{
\begin{minipage}{5cm}{
\centering
  \begin{pspicture}(-1.5,-1.5)(1.5,1.5)
    \psarc(0,0){1}{0}{90}
    \psarc(0,0){1}{180}{270}
    \psarc[linestyle=dotted](0,0){1}{90}{180}
    \psarc[linestyle=dotted](0,0){1}{270}{360}
    \rput(1.2,.7){$\partial_D \ball$}
    \rput(-1.2,.7){$\partial_N \ball$}
    \rput(-1.2,-.7){$\partial_D \ball$}
    \rput(1.2,-.7){$\partial_N \ball$}
  \end{pspicture}}
\end{minipage}
\begin{minipage}{5cm}{
\centering
  \begin{pspicture}(-1.5,-1.5)(1.5,1.5)
    \psarc(0,0){1}{0}{60}
    \psarc(0,0){1}{120}{180}
    \psarc(0,0){1}{240}{300}
    \psarc[linestyle=dotted](0,0){1}{60}{120}
    \psarc[linestyle=dotted](0,0){1}{180}{240}
    \psarc[linestyle=dotted](0,0){1}{300}{360}
    \rput[l](.9,.6){$\partial_D \ball$}
    \rput[r](-.9,.6){$\partial_D \ball$}
    \rput(0,-1.2){$\partial_D \ball$}
    \rput[r](-.9,-.6){$\partial_N \ball$}
    \rput[l](.9,-.6){$\partial_N \ball$}
    \rput(0,1.2){$\partial_N \ball$}
  \end{pspicture}
}
\end{minipage}}
\end{minipage}
\caption{Uniform $2$--partition (left) and uniform $3$--partition (right).} \label{unifKpart}
\end{figure}

\subsection{The first mixed eigenvalue}

The extrema of the first mixed eigenvalue has been studied by Denzler
\cite{denzler:bounds}. In a later paper \cite{denzler:disc} he studied the same problem with a particular focus on minimization and disks of any
dimension. Denzler showed that the minimal first eigenvalue of the $d$--dimensional ball is achieved when the Dirichlet's part $\partial_D \Omega$ is a spherical cap. In dimension $2$, this is just one connected part, that is, the uniform $1$--partition. Cox and Uhlig give another proof of this result in \cite{cu:drum}.

Burchard and Denzler investigated the case of the first mixed
eigenvalue of $2$--dimensional domains, particularly the square
\cite{bd:square}. In this case, it seems that there is no absolute minimizing arrangement at all; minimizer depends on the family of boundary conditions where the minimizer is taken.

It is known that the maximizing arrangement of boundary
conditions is not attained and that one should ``smear'' the boundary
condition in order to increase the eigenvalues. In Cox and Uhlig \cite{cu:drum} and Denzler \cite{denzler:bounds} this result is shown for the first eigenvalue of the mixed Laplacian on a $d$--dimensional Lipschitz domain. As we shall see in Section~\ref{sectionMAX} this result can be extended to any higher eigenvalues (Theorem \ref{thm:maxk}).

Moreover, we prove that on the disk, with the added constraint that number of connected components of the Dirichlet part is bounded above by $n$, then the maximum of the first eigenvalue is attained when these parts are uniformly distributed in the border (Theorem \ref{thm:max1:k-com}).

\begin{theorem}
  \label{thm:max1:k-com}
  The uniform $n$--partition of length $\ell$ is a maximizer for the first
  mixed eigenvalue on the disk among all parts of boundary of length $\ell$ with
  at most $n$ connected components.
\end{theorem}

\subsection{Minimizing the second mixed eigenvalue of the disk.} Let us state the following Conjecture which will be justified in this paper.

\begin{conjecture}
  \label{conj:min2}
  The minimizing arrangement of boundary conditions for the second eigenvalue of the mixed Dirichlet--Neumann problem on the disk is given by the uniform $2$--partition.
\end{conjecture}

Numerics in Section~\ref{Morenum} support this conjecture and we are able to prove some partial results which indicate that it should be true. The following theorem reduces the problem of finding a minimum among a smaller family of boundary conditions containing the uniform $2$--partition.

\begin{theorem}
  \label{thmB}
  The minimal second mixed eigenvalue of the disk, among all Dirichlet parts of length $\ell$, is attained on some $\Gamma^* \in \mF_\ell$, which has at most two connected components. Moreover, if $\Gamma^*$ has exactly two connected components, they are of equal length.
\end{theorem}

That is, the minimal second mixed eigenvalue lies in the following family.

\begin{definition}[Boundary family $\mF_\ell^e$]
  \label{def:family-e}
 Let $\mF_\ell^e$ be the family of Dirichlet's parts which have at most two connected components and these components. If there are two then they are of equal length.
\end{definition}

 Such partitions are necessarily symmetric with respect to a single line of symmetry and the uniform 2--partition is the only one of these which is symmetric with respect to two (perpendicular) axes of symmetry. As we shall see in Section \ref{seclem}, for a given $\ell$ this family can be parameterized by a compact interval. Theorem \ref{thmB} claims that the minimizer of the second eigenvalue among $\mF_{\ell}$ lies in $\mF_{\ell}^e$. Note that $\Gamma_2 \in \mF_{\ell}^e$. More precisely, we prove in Section \ref{seclem} (Lemma \ref{lemSym}) that the second mixed eigenvalue of any boundary condition on the disk is higher than one of the second mixed eigenvalue associated to this family.

 It seems that one of the difficulties in obtaining a more precise result, that is to prove Conjecture 1, is the lack of knowledge of the behavior of nodal lines in the case of mixed boundary conditions. As we shall see, the behavior of the nodal line of an eigenfunction associated to the second eigenvalue, usually called the first nodal line, is related to the minimizer in some sense, see Remark~\ref{remLONG}. But in the cases of mixed boundary condition we do not even know if the nodal line is closed or if it reaches the border, that is Payne's Conjecture for the mixed cases is not yet proved, see Theorem~\ref{thm:nodal} for a partial result.

 Moreover, usual techniques and tools of optimization using derivatives cannot be used directly in the case of mixed boundary conditions because mixed eigenfunctions may be non-smooth on the border of the domain, see~\cite{miranda:mbvp}. In particular, we can not extend directly Melas's or Alessendrini's proof of the Payne's Conjecture for the mixed cases.

 Another difficulty is the rigidity of the problem, that is we must work on a fixed domain and with functions that we do not know explicitly. The usual technique of finding minimums of low eigenvalues is firstly to do a rearrangement (e.g. a spherical symmetrization) of eigenfunctions which decreases their Rayleigh quotient and then to make use of the variational characterization of eigenvalues. However, there are not many rearrangements that one can use here since one has to obtain functions defined on the disk and which satisfy some orthogonality relation.

{\bf Acknowledgements} -- The author is very grateful to her M.Sc. supervisor, Professor Iosif Polterovich, for comments, helpful discussions and suggesting this problem.

The author wishes to thank Professor Michael Levitin for valuable comments, as well as his invitation and hospitality at Heriot--Watt University. The author thanks Michael K. Graham and also Claude Gravel for his help with MatLab.
\section{Minimizing the second mixed eigenvalue of the disk}

\subsection{A key lemma}
\label{seclem}

We now assume that the domain is the unit disk, $\ball \subset \R^2$, and we denote
\begin{equation}
  \lambda_k(\partial_D \ball) := \lambda_k(\ball, \partial_D \ball).
\end{equation}

The next lemma implies the existence of a minimizer and reduces the
problem of its location to the family $\mF_\ell^e$. First, we give another description of this family, see Definition~\ref{def:family-e}.

\begin{figure}
  \centering
    \begin{pspicture}(-1,-1.5)(1.5,1.5)
      \psarc[linestyle=dotted](0,0){1}{-30}{30}
      \psarc(0,0){1}{30}{130}
      \psarc[linestyle=dotted](0,0){1}{130}{230}
      \psarc(0,0){1}{230}{330}
      \psline[linestyle=dotted](0.866,0.5)(0,0)
      \psline[linestyle=dotted](1,0)(0,0)
      \rput(.7,.2){$\beta$}
    \end{pspicture}
    \caption{Boundary condition $\Gamma(\beta)$ }
    \label{fig:Gamma}
\end{figure}
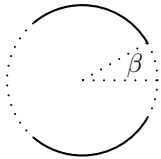

\begin{definition} \label{betafamily} $\mF_\ell^e = \{ \Gamma(\beta) | \, \beta \in [0,\frac{2\pi-\ell}{4}]\}$, where

\begin{align}
  \Gamma(\beta) := \left\{
    e^{i\theta} \left| \,
    \begin{array}{l}
    \theta \in [0,2\pi],\; {\rm and}\\
    \theta \notin \left(-\beta,\beta \right)
    \cup \left(\pi-\left( \frac{\ell}{2}+\beta \right), \pi + \left( \frac{\ell}{2}+\beta \right) \right) \\
    \end{array}
    \right.\right\}
\end{align}
(see Figure~\ref{fig:Gamma}).
\end{definition}

Note that the uniform $2$--partition is given by $\Gamma(\frac{2\pi-\ell}{4}) = \Gamma_2$ and the uniform $1$--partition is given by $\Gamma(0)=\Gamma_1$.

Recall that $\ball$ is the {\it open} unit disk centered at the origin of $\bR^2$. For a mixed problem on $\ball$, determined by a Dirichlet part $\partial_D \ball \in \mF_{\ell}$, and an associated second eigenfunction $u_2$, the nodal domains are the connected, open sets $$\ball_+(u_2) := \{ p \in \ball | u_2(p) > 0 \} \;\;\; \text{and} \;\;\; \ball_-(u_2) := \{ p \in \ball | u_2(p) < 0 \}.$$

\begin{definition}
 Given $\partial_D \ball \in \mF_{\ell}$ and a corresponding second
  eigenfunction $u_2$ with nodal domains $\ball_+$ and $\ball_-$, we define
  \begin{equation}
    \beta_{u_2} := \min\left\{
      \left|\partial_N \ball \cap \partial\ball_+\right|,
      \left|\partial_N \ball \cap \partial\ball_-\right|
    \right\}/2.
  \end{equation}
\end{definition}

\begin{lemma}
  \label{lemSym} Given $\partial_D \ball \in \mF_{\ell}$ and any associated second eigenfunction $u_2$, we have $\lambda_2(\partial_D \ball) \geq \lambda_2(\Gamma(\beta_{u_2}))$, with equality if and only if $\partial_D \ball = \Gamma (\beta_{u_2})$ up to a rotation,
  a reflection and a set of zero measure.
\end{lemma}

\begin{rem}
\label{remLONG}
There are two extremal cases worth noting :

 \begin{enumerate}
  \item If the nodal line of $u_2$ divides $\partial_D \ball$ in two equal parts then $\beta(u)=\frac{2\pi-\ell}{4}$. Hence Lemma~\ref{lemSym} gives that $\lambda_2(\partial_D \ball) \geq \lambda_2(\Gamma_2)$ where $\Gamma_2$ stands for the uniform $2$--partition.

  \item If the nodal line of $u_2$ intersects the boundary of $\ball$ in at most one point, then $\beta(u)=0$ and by Lemma~\ref{lemSym}, $\lambda_2(\partial_D \ball) \geq \lambda_2(\Gamma_1)$, where $\Gamma_1$ stands for the uniform $1$--partition.
 \end{enumerate}
\end{rem}
\begin{proof}
  Define the restrictions of $u_2$ to $\ball_+$ and $\ball_-$ as $u_+$ and
  $u_-$ respectively. Taking $-u_2$ if necessary, we can suppose that $2\beta_{u_2} = \left|\partial_N \ball \cap \partial\ball_+\right| \leq \left|\partial_N \ball \cap
    \partial\ball_-\right|$. We denote the nodal
  line $$\mathcal{N} := \{(x,y) \in \ball | u_2(x,y) = 0 \} = \ball \backslash (\ball_+ \cup \ball_- ).$$

  The function $u_+$ is a first eigenfunction on the domain $B_+$ with Dirichlet
  boundary condition on $\mathcal{N} \cup \partial_D \ball \cap \partial
  \ball_+$ and Neumann boundary condition on the remainder. Now use spherical symmetrization, as described by Polya and Szegõ in~\cite{ps:sym}, to
  rearrange the function $u_+$ in the angular direction with respect to the
  positive $x$--axis and centered at the origin. The result is a smooth function $u_+^*$ defined on the domain $\ball_+^*$. Recall, that the spherical symmetrization with respect to the origin and the positive $x$--axis implies that for any circle $S_r$ of radius $r \leq 1$ centered at the origin, we have $\left|S_r \cap \ball_+\right| =\left|S_r \cap \ball_+^*\right|$ and $S_r \cap \ball_+^*$ is connected, symmetric with respect to the $x$--axis and passes through the positive $x$--axis (see Figures~\ref{fig:ballp} and~\ref{fig:ballpstar}~. Moreover, the function $u_+^*$ is symmetric with respect to the $x$--axis and, on $S_1$, it is a strictly positive function for $\theta \in (-\beta_{u_2}, \beta_{u_2})$ and vanishes elsewhere on the boundary of $\ball_+^*$.

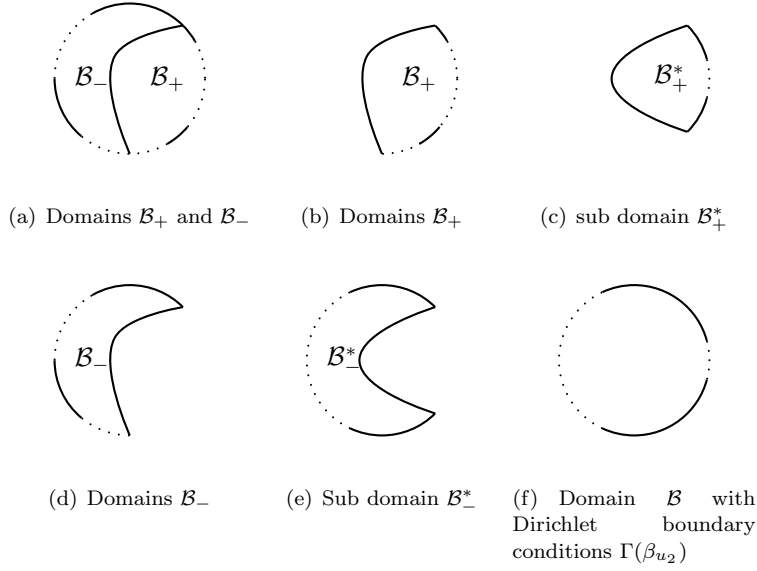
\begin{figure}[ht]
  \centering
  \subfigure[Domains $\ball_+$ and $\ball_-$]{
    \label{fig:ballpm}
    \begin{pspicture}(-1.5,-1.5)(1.5,1.5)
      \psarc[linestyle=dotted](0,0){1}{0}{30}
      \psarc(0,0){1}{45}{120}
      \psarc(0,0){1}{180}{230}
      \psarc(0,0){1}{300}{320}
      \psarc(0,0){1}{30}{45}
      \psarc[linestyle=dotted](0,0){1}{120}{180}
      \psarc[linestyle=dotted](0,0){1}{230}{300}
      \psarc[linestyle=dotted](0,0){1}{320}{360}
      \pscurve(0,-1)(-0.2,0.3)(.707,.707)
      \rput(.5,0){$\ball_+$}
      \rput(-.5,0){$\ball_-$}
    \end{pspicture}
  }
      \subfigure[Domains $\ball_+$]{
    \label{fig:ballp}
    \begin{pspicture}(-1.5,-1.5)(1.5,1.5)
      \psarc[linestyle=dotted](0,0){1}{0}{30}
      \psarc(0,0){1}{300}{320}
      \psarc(0,0){1}{30}{45}
      \psarc[linestyle=dotted](0,0){1}{270}{360}
      \pscurve(0,-1)(-0.2,0.3)(.707,.707)
      \rput(.5,0){$\ball_+$}
    \end{pspicture}
  }
  \subfigure[sub domain $\ball_+^*$]{
    \label{fig:ballpstar}
    \begin{pspicture}(-1.5,-1.5)(1.5,1.5)
      \psarc[linestyle=dotted](0,0){1}{-15}{15}
      \psarc(0,0){1}{15}{45}
      \psarc(0,0){1}{315}{345}
      \pscurve(0.707,-0.707)(-.3,0)(.707,.707)
      \rput(.5,0){$\ball_+^*$}
    \end{pspicture}
  }
    \subfigure[Domains $\ball_-$]{
    \label{fig:ballm}
    \begin{pspicture}(-1.5,-1.5)(1.5,1.5)
      \psarc(0,0){1}{45}{120}
      \psarc(0,0){1}{180}{230}
      \psarc[linestyle=dotted](0,0){1}{120}{180}
      \psarc[linestyle=dotted](0,0){1}{230}{270}
      \pscurve(0,-1)(-0.2,0.3)(.707,.707)
      \rput(-.5,0){$\ball_-$}
    \end{pspicture}
  }
    \subfigure[Sub domain $\ball_-^*$]{
    \label{fig:ballmstar}
    \begin{pspicture}(-1.5,-1.5)(1.5,1.5)
      \psarc(0,0){1}{45}{115}
      \psarc[linestyle=dotted](0,0){1}{115}{245}
      \psarc(0,0){1}{245}{315}
      \pscurve(0.707,-0.707)(-.3,0)(.707,.707)
      \rput(-.5,0){$\ball_-^*$}
    \end{pspicture}
  }
  \subfigure[Domain $\ball$ with Dirichlet boundary conditions
  $\Gamma(\beta_{u_2})$]{
    \label{fig:ballgamma}
    \begin{pspicture}(-1.5,-1.5)(1.5,1.5)
      \psarc[linestyle=dotted](0,0){1}{-15}{15}
      \psarc(0,0){1}{15}{115}
      \psarc[linestyle=dotted](0,0){1}{115}{245}
      \psarc(0,0){1}{245}{345}
    \end{pspicture}
  }
  \caption{Domains used in proof of Lemma~\ref{lemSym}}
  \label{fig:balls}
\end{figure}

 Similarly, $u_-$ is a first eigenfunction on the domain $B_-$ with Dirichlet boundary condition on $\mathcal{N} \cup \partial_D \ball \cap \partial \ball_-$ and Neumann boundary condition on the remainder. We proceed in a similar way to rearrange $u_-$ but the spherical symmetrization should be done with respect to the negative $x$--axis. So we get the function $u^*_-$ defined and smooth on $\ball^*_-$, strictly negative on an arc in $S_1$ of length $|\del_N \ball \cap \ball_-| = 2\pi -\ell -2\beta_{u_2}$ and vanishing elsewhere on $\del \ball_-$, see Figures~\ref{fig:ballm} and~\ref{fig:ballmstar}.

  Since we have $\left|S_r \cap \ball_+\right|
  =\left|S_r \cap \ball_+^*\right|$ and $\left|S_r \cap \ball_-\right|
  =\left|S_r \cap \ball_-^*\right|$ for any circle, $S_r$, the interior of
  $\overline{\ball_+^*} \cup \overline{\ball_-^*}$ is a disk of radius one. Remark also that $\{p \in S_1 | u_+(p)=0 \text{ or } u_-(p)=0 \} = \Gamma(\beta_{u_2})$.

 By a result due to Sperner~\cite{sp:sym}, the spherical symmetrization process ensures that
  \begin{equation}
    \begin{split}
      \label{eq:greatfirstplus}
      \lambda_2(\partial_D \ball) &= \lambda_1(\ball_+, \mathcal{N} \cup (\del_D \ball \cap \partial \ball_+))= \frac{\int_{\ball_+} \abs{ \grad
          u_+}^2}{\int_{\ball_+} \abs{u_+}^2} \\
      &\geq \frac{\int_{\ball_+^*}
        \abs{ \grad u_+^*}^2}{\int_{\ball_+^*} \abs{u_+^*}^2} \geq \lambda_1(\ball_+^*, \mathcal{N}^* \cup (\Gamma(\beta_{u_2}) \cap
      \partial \ball_+^*))
    \end{split}
  \end{equation}
  and similarly for the domain $\ball_-^*$ we have
  \begin{equation}
    \begin{split}
      \label{eq:greatfirstminus}
      \lambda_2(\partial_D \ball) &= \lambda_1(\ball_-, \mathcal{N} \cup (\del_D \ball \cap \partial \ball_-))= \frac{\int_{\ball_-} \abs{ \grad
          u_-}^2}{\int_{\ball_-} \abs{u_-}^2} \\
      &\geq \frac{\int_{\ball_-^*}
        \abs{ \grad u_-^*}^2}{\int_{\ball_-^*} \abs{u_-^*}^2} \geq \lambda_1(\ball_-^*, \mathcal{N}^* \cup (\Gamma(\beta_{u_2}) \cap
      \partial \ball_-^*)).
    \end{split}
  \end{equation}

But the Theorem of Domain Monotonicity of Eigenvalues, for Dirichlet data (see e.g.~\cite{chavel:eigen}) implies that

  \begin{equation}
    \begin{split}
      \label{eq:lessmax}
      \lambda_2(\Gamma(\beta_{u_2})) \leq \max\bigl\{&
      \lambda_1(\ball_+^*, \mathcal{N}^* \cup (\Gamma(\beta_{u_2}) \cap
      \partial \ball_+^*) ,\\
      &\lambda_1(\ball_-^*, \mathcal{N}^* \cup
      (\Gamma(\beta_{u_2}) \cap \partial \ball_-^*)) \bigr\}.
    \end{split}
  \end{equation}

Hence, Inequalities~\eqref{eq:lessmax},
  \eqref{eq:greatfirstplus}, and \eqref{eq:greatfirstminus} complete the
  proof that $\lambda_2(\Gamma(\beta_{u_2})) \leq \lambda_2(\partial_D
  \ball)$.

~\\
{\it The case of equality} -- The equality $\lambda_2(\Gamma(\beta_{u_2})) = \lambda_2(\partial_D \ball)$ implies equality in~\eqref{eq:lessmax}. Denote $\partial_D \ball^*_{\pm} := \mathcal{N}^* \cup (\Gamma(\beta_{u_2}) \cap \partial \ball_{\pm}^* $ and $\partial_N \ball^*_{\pm}$ the respective reminder.

We first prove that, in this case, we must have $$\lambda_1(\ball_+^*, \partial_D \ball^*_+) = \lambda_1(\ball_-^*,\partial_D \ball^*_-).$$

  Indeed if it is not, we can suppose without loss of generality that $$\lambda_2(\Gamma(\beta_{u_2})) = \lambda_1(\ball_+^*, \partial_D \ball^*_+) > \lambda_1(\ball_-^*, \partial_D \ball^*_-).$$ Thus, by the continuity of the first mixed eigenvalue with respect to continuous variation of the Dirichlet part of the boundary, there exists a domain $D \subset \ball$ such that $\ball_+^* \subset D $ and $ \ball \backslash \overline{D} \subset \ball_-$ in such a way that $\lambda_1(D,\del D \backslash \partial_N \ball^*_+ ) = \lambda_1(\ball \backslash \overline{D}, \del (\ball \backslash \overline{D}) \backslash \partial_N \ball^*_- )$. On another hand, Monotonicity of Eigenvalues and Domain Monotonicity of Eigenvalues (for Dirichlet Data) imply that \begin{equation}
    \begin{split}\lambda_2(\Gamma(\beta_{u_2})) &= \lambda_1(\ball_+^*, \partial_D \ball^*_+)\\ &> \lambda_1(D,\del D \backslash \partial_N \ball^*_+ )\\ &= \lambda_1(\ball \backslash \overline{D}, \del (\ball \backslash \overline{D}) \backslash \partial_N \ball^*_- )\\ &\geq \lambda_2(\Gamma(\beta_{u_2})).    \end{split}
  \end{equation} This is a contradiction.

  This proves that the equality $\lambda_2(\Gamma(\beta_{u_2})) = \lambda_2(\partial_D \ball)$ implies equality in {\it both} \eqref{eq:greatfirstplus} and \eqref{eq:greatfirstminus}. But this can occur if and only if $\partial_N \ball_{\pm}$ equals $(\partial_N \ball_{\pm})^*$ (and $\ball_{\pm}$ equals $\ball_{\pm}^*$) up to rotation, reflection and set of zero measure. Indeed, this can be deduced directly (in view of the definition of $u_{\pm}$) from the proof given by Denzler in~\cite{denzler:disc} concerning the case of equality. Therefore, we have $\partial_D\ball = \Gamma(\beta_{u_2})$ up to a rotation and a set of zero measure. \end{proof}

\subsection{Corollaries and Extension of Lemma \ref{lemSym}}
\label{sectTHEO}
~\\

{\bf Extension to higher eigenvalues} -- The proof of Lemma~\ref{lemSym} uses rearrangements. A natural question is wether or not the proof can be adapted to higher eigenvalues. Since we use the fact that the number of nodal domains of a second eigenfunction is exactly 2, the adaptation does not seem straightforward.\\

{\bf Extension to higher dimensions} -- Lemma~\ref{lemSym} could be extended to higher dimensions. In this case the minimizing arrangement with given volume must have a Neumann part consisting of two opposite spherical caps with possibly different volumes.

More precisely, let $\Gamma_{\ell_N}(\beta) \subset S^{n-1}$ be the union of two opposite spherical caps of volume $\beta$ and $\ell_N - \beta$ (where $\beta < \ell_N/2$ and $\ell_N \in [0,\text{Vol}(S^{n-1})]$). Consider a mixed eigenvalue problem on the ball $\ball^n$ in $\bR^n$ with Dirichlet boundary condition on $\del_D \ball^n \subset S^{n-1}$ and Neumann condition on the open reminder $\del_N \ball^n = S^{n-1} \backslash \overline{\del_D \ball^n}$ with $|\text{Vol}(\del_N \ball^n)|=\ell_N$. A second eigenfunction, $u_2$, associated to this problem has exactly two nodal domains $\ball^n_+$ and $\ball^n_-$. As above, let

\begin{equation}
    \beta_{u_2} := \min\left\{
      \left|\partial_N \ball^n \cap \partial\ball^n_+\right|,
      \left|\partial_N \ball^n \cap \partial\ball^n_-\right|
    \right\}/2.
  \end{equation}

The proof of Lemma~\ref{lemSym} works with spherical symmetrization in this framework as well.

\begin{lemma}[extension of Lemma~\ref{lemSym}] Given $\partial_D \ball^{n} \in \mF_{\ell}$ and any associated second eigenfunction $u_2$, we have $$\lambda_2(\ball^{n},\del_D \ball^n) \geq \lambda_2(\ball^{n},S^{n-1} \backslash \overline{\Gamma_{\ell_N}(\beta_{u_2})}),$$ with equality if and only if $\partial_D \ball^{n} = \Gamma (\beta_{u_2})$ and $\del_N \ball^n = \Gamma_{\ell_N}(\beta_{u_2})$ up to isometry and set of zero measure.
\end{lemma}

{\bf Theorem \ref{thmB}} -- Theorem \ref{thmB} is a straightforward corollary of Lemma~\ref{lemSym}. Recall that it gives geometric properties of the minimizer of the second eigenvalue among all Dirichlet's parts of fixed length.

\begin{proof}[Proof of Theorem  \ref{thmB}] Lemma~\ref{lemSym} gives that for each $\partial_D \ball \in \mF_{\ell}$, there exists $\Gamma(\beta) \in \mF_\ell^e$ such that $\lambda_2(\partial_D \ball) \geq \lambda_2(\Gamma (\beta))$ with equality if and only if $\partial_D \ball = \Gamma (\beta)$ up to a set of zero measure and an isometry. The inequality implies that the minimum among this family is a global minimum and the case of equality implies that if there is a minimizer, it lies in this family.
\end{proof}

\begin{corollary}[existence]
There exists a minimizer of the second eigenvalue among all Dirichlet parts of fixed length.
\end{corollary}

\begin{proof}
One can show that the map $\beta \mapsto \lambda_2(\Gamma (\beta))$ is continuous via standard methods, such as monotonicity of mixed eigenvalues and continuity of the first mixed eigenvalue. Hence, the compactness of the interval of definition of $\beta$ implies the existence of a minimum in this family and Lemma~\ref{lemSym} implies that this is a minimum among all Dirichlet parts of fixed length.
\end{proof}

{\bf Nodal line and the family $\mF_\ell^e$ } -- The rearrangement used in the proof of Lemma~\ref{lemSym} shows that it could be useful to understand the behavior of nodal lines. For example, if the first nodal line of any mixed problem associated to a Dirichlet's part in $\mF_\ell^e$ is closed then it follows from Lemma~\ref{lemSym} that Conjecture \ref{conj:min2} is false (see the second case in Remark~\ref{remLONG}). On the other hand, Lemma~\ref{lemSym} implies that the second eigenvalue associated to the uniform $2$--partition is lower than any second eigenvalue associated to mixed problem on the disk for which the first nodal line divides equally the Neumann's part of the boundary (see the first case in Remark~\ref{remLONG}). It is, therefore, a natural question to ask wether or not the first nodal lines are closed.

The next theorem claims that the answer is no: for any mixed problem on the disk associated to a Dirichlet's part in $\mF_\ell^e$ the associated first nodal line cannot be closed. This theorem is a particular case of a result announced in~\cite{Gr:th} which proposed a mixed counterpart to Lin's result~\cite{lin:conj} about Payne's Conjecture in the case of pure Dirichlet problem for symmetric, smooth and convex domains. The proof below is given for the paper to be self-contained and follows the idea appearing in the proof of Payne~\cite{payne:conj}.

\begin{theorem}
  \label{thm:nodal}
Consider a mixed boundary value problem on a disk such that the Dirichlet and Neumann parts of the boundary
$\Gamma_D$ are axially symmetric. Let $u$ be an eigenfunction corresponding to the
second eigenvalue $\lambda_2$. Then the nodal line $\mN(u) = \overline{ \{x \in \Omega \; | \; u(x) = 0\}}$ of u cannot be closed.
\end{theorem}
\begin{proof} Let $\Gamma_D$ be symmetric with respect to the $x$-axis, denote 
$$ \Gamma_D^+ := \{(x,y) \in \Gamma_D \;|\; y > 0\} \;\;\;  \text{and} \;\;\; \Gamma_D^- := \{(x,y) \in \Gamma_D \;|\; y < 0\}$$ and let $u$ be a second mixed eigenfunction which have a closed nodal line $\mN(u)$ (which thus does not meet the boundary). Thus $u$ may be chosen symmetric with respect to the $x$--axis and the function $v:= \frac{\del u}{\del y}$ is well-defined and anti-symmetric on $\ball \cup \Gamma_D$. Define $\ball_*^+ := \{ (x,y) \in \ball \;|\; v(x,y) < 0 \;\; \text{and} \;\; y>0 \}$
and $ \ball_* := \{ (x,y) \in \ball \;|\; (x,\pm y) \in \ball_*^+ \}$. The restriction of the function $v$ on $\ball_*$ is such that :
\begin{itemize}
\item $-\Delta v = \lambda_2(\Gamma_D) v$ (since the Laplacian commutes with the $y$--derivative)
\item $v$ is $L^2$--orthogonal to the first eigenfunction of any mixed problem on $\ball_*$ (because it is anti-symmetric)
\item $v=0$ on $\del \ball_* \cap \ball$
\end{itemize}

Hence $v$ is a test function for the second mixed eigenvalue $\lambda_2(\ball_*, \del \ball_* \cap \ball)$ and
\begin{equation}
 \label{eq:lamb+}
\lambda_2(\Gamma_D) \geq \lambda_2(\ball_*, \del \ball_* \cap \ball).
\end{equation}

We want to prove that $\ball_*$ is non-empty and does not contain an open part of $\Gamma_D$ in its closure, in order to use Monotonicity of mixed eigenvalue.

By hypothesis the first nodal line $\mN(u)$ is closed, so:
\begin{itemize}
\item $\mN(u)$ divides $\ball$ in $2$ disjoint connected components : the inside part, whose boundary is only $\mN(u)$ and the outside part, whose boundary contains also the boundary of the disk
\item we can assume that $u$ is positive on the inside part and negative on the outside one
\item $\mN(u) \cap \del \ball = \emptyset$, so $u$ is smooth around $\mN(u)$ and $\mN(u)$ is a $C^2$-immersed circle without intersection, see~\cite{sy:lect}.
\end{itemize}
Denote by $\frac{\del}{\del \eta}$ the outward unit normal vector field to the boundary $\mN(u)$ of the inside part and by $\frac{\del}{\del n}$ the outward unit normal vector field to the boundary $\del \ball$ of the disk. Then $\frac{\del u}{\del \eta} < 0$ on $\mN(u)$ and $\frac{\del u}{\del n} > 0$ on $\Gamma_D$ (Hopf Boundary Lemma).

Hence the set $\ball_*$ is not empty since there must exist a point of $\mN(u)$, in the upper half plane, where $\frac{\del }{\del y}$ and $\frac{\del }{\del \eta}$ are collinear and have the same orientation.

Denote by $ \frac{\del }{\del \theta}$ a vector field, tangent to $\Gamma_D$, so there exist $a$ and $b \in C^{\infty}(\Gamma_D)$ such that $\frac{\del }{\del y} = a\frac{\del }{\del n} + b\frac{\del }{\del \theta}$ on $\Gamma_D$. Since $u$ is constant on $\Gamma_D$, for any $P \in \Gamma_D$  $$v(P)=\frac{\del u}{\del y} (P)= a(P)\frac{\del u}{\del n} (P) + b(P)\frac{\del u}{\del \theta} (P) = a(P)\frac{\del u}{\del n}(P).$$

Observe moreover that $ a= \langle\frac{\del }{\del y}, \frac{\del }{\del n} \rangle  > 0$ on $\Gamma_D^+$ and $ a= \langle\frac{\del }{\del y}, \frac{\del }{\del n} \rangle < 0$ on $\Gamma_D^-$ and as stated before $\frac{\del u}{\del n} > 0$ on $\Gamma_D$. Hence
\begin{equation}
v(P) = \left\{
\begin{array}{l}
>0  \;\;\; \text{for} \; P \in \Gamma_D^+, \\
<0  \;\;\; \text{for} \; P \in \Gamma_D^-.
\end{array}
\right.
\end{equation}

 Thus  $\del \ball_* \cap \Gamma_D = \emptyset$ or is only a finite set (the intersection of the $x$-axis and $\Gamma_D$). In particular, $\ball_*$ is not $\ball$ itself.

  So we can use Monotonicity (with Dirichlet data~\cite{chavel:eigen}) and we get $$\lambda_2(\ball_*, \del \ball_* \cap \ball) > \lambda_2(\Gamma_D)$$ which, combined with inequality~\eqref{eq:lamb+}, is a contradiction.
\end{proof}

\begin{rem} The last proof can easily be extended to more general situations and has been presented in this case to simplify the reading. We stated Theorem~\ref{thm:nodal} with the hypothesis that the domain is the disk and $\Gamma_D$ is symmetric with respect to one axis, but the hypothesis effectively used are similar to Payne's ones: the domain is symmetric with respect to the $x$--axis and is $y$--convex (i.e $\langle \frac{\del}{\del y}, \frac{\del}{\del n}\rangle$ has a constant sign on each side of the $x$--axis, see~\cite{payne:conj}) and $\Gamma_D$ is symmetric with respect to $x$--axis, see also~\cite{Gr:th}.
\end{rem}

\section{Numerics}\label{Morenum}
Theorem \ref{thmB} states that among all partitions of the boundary of given length,
the minimizer of the second eigenvalue of the disk belongs to the family
$\mF^e_{\ell}$ (see Definition~\ref{def:family-e} and \ref{fig:Gamma}). Hence, we now focus on this particular family.

\begin{figure}[h]
\centering
\begin{minipage}{4.25in}{

\begin{minipage}{4.2in}{\includegraphics{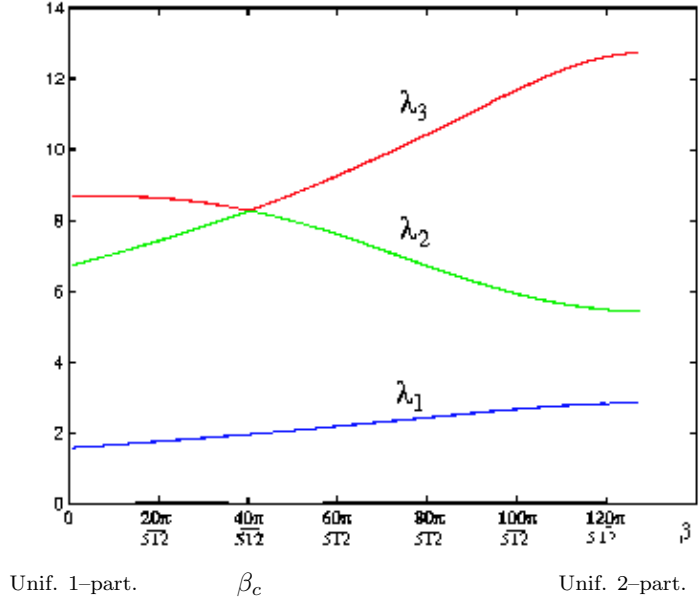}}\end{minipage}\\
\begin{minipage}{0.5cm}{$\,$}\end{minipage}
\begin{minipage}{1in}{\footnotesize{Unif. $1$--part.}}\end{minipage}
\begin{minipage}{0.25cm}{$\,$}\end{minipage}
\begin{minipage}{1in}{$\beta_c$}\end{minipage}
\begin{minipage}{1.5cm}{$\,$}\end{minipage}
\begin{minipage}{1in}{\footnotesize{Unif. $2$--part.}}\end{minipage}
}
\end{minipage}
\caption{Variation of first three eigenvalues in function to $\beta$ for $\ell = \pi$.}\label{graph3eig}
\end{figure}

Figure \ref{graph3eig} shows the graph of the first three eigenvalues over the family $\mF^e_{\pi}$, computed for $\beta = \{\frac{i\pi}{512}\;|\; i=0, 1, \dots, 128\}$. The graph indicates that, in this family, the minimum of the second eigenvalue is reached for $\beta=\frac{2\pi-\ell}{4}= \frac{\pi}{4}$. Since the corresponding partition is $\Gamma(\frac{2\pi-\ell}{4}) =\Gamma_2$, the uniform $2$--partition of length $\ell$. This numerical test supports Conjecture~1.

Let us denote by $\beta_c$, the point in the family where the two curves representing $\lambda_2$ and $\lambda_3$ on Figure \ref{graph3eig} seems to intersect. It is the only element of the family for which the multiplicity of $\lambda_2$ could be two. Moreover, this critical point corresponds to an important change in the position of the corresponding nodal line. In fact, for $0 \leq \beta < \beta_c$, simulations show that the nodal line is horizontal with respect to the orientation given in Figure \ref{fig:Gamma} and then, the corresponding eigenfunction is anti-symmetric. For $\beta_c < \beta \leq \frac{\pi}{2}$ the nodal line is vertical according to numerics and the second eigenfunction should be symmetric with respect to the line of symmetry of the Family $\mF^e_{\pi}$.

\begin{figure}[h]
\psfrag{b}[][][1]{$\beta = \frac{350\pi}{1024}$}
\psfrag{a}[][][1]{$\beta = \frac{351\pi}{1024}$}
\includegraphics{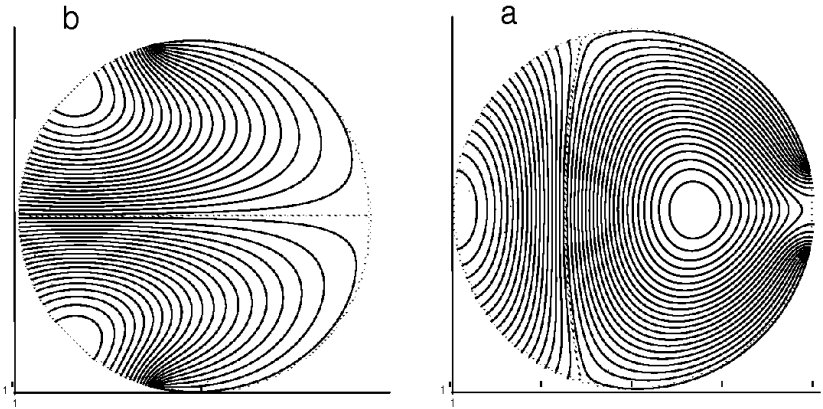}
\caption{Level lines of a second eigenfunction for $\beta$ near $\beta_c$.}
\label{nodalline}
\end{figure}

Figure ~\ref{nodalline} shows level lines of only two tests, namely for $b = \frac{350\pi}{1024} < \beta_c$ and $a = \frac{351\pi}{1024} > \beta_c$. The nodal line is represented by the doted line. This suggest that before $\beta_c$ second eigenfunctions are anti-symmetric and after $\beta_c$ second eigenfunctions are symmetric. Then, Lemma \ref{lemSym}, the first case of Remark~\ref{remLONG} and previous numerical evidence seem to indicate that the main difficulty in establishing a proof of Conjecture 1 is showing the result for $\beta \geq \beta_c$.

\section{Maximization}\label{sectionMAX}

We now study the existence and geometric properties of an arrangement of boundary condition maximizing some mixed eigenvalue.

 In~\cite{cu:drum} and~\cite{denzler:bounds}, it is shown that the first eigenvalue of the mixed Laplacian has no maximizing arrangement on $d$--dimensional Lipschitz domain and that the first Dirichlet eigenvalue is the supremum. The following theorem is a natural extension of this result for any eigenvalue and our proof is a straightforward adaptation of Denzler's.

  \begin{theorem}
  \label{thm:maxk}
  Let $\Omega$ be a bounded Lipschtiz domain with boundary $\partial \Omega$, and fix $ 0 < A < |\partial\Omega|$. Then
  \begin{equation}
    \sup\{\lambda_k(\Omega,
    \partial_D \ball)\; | \; |\partial_D \Omega|= A\}=
    \lambda_k(\Omega,\partial\Omega)
  \end{equation}
  and a maximizing sequence of parts can be given in such a way that their
  characteristic functions converge weakly to a constant in
  $L^2(\partial\Omega)$.
\end{theorem}

  Hence, in order to increase any mixed eigenvalue one should {\it smear} the boundary condition as much as possible. This concludes the question of maximization of mixed eigenvalue without other constraints than fixed volume. Nevertheless, we prove below that uniform $n$--partition maximizes the first eigenvalue under the constraint that number of connected components is bounded above by $n$. It seems to be a natural specification of the Theorem~\ref{thm:maxk} for the case $k=1$ and on the disk. In our notation and for the disk, such a {\it smearing} takes the form of the uniform $n$--partition and we would expect the first eigenvalue to increase with $n$. However as the Conjecture 1, Numerics of Section~\ref{Morenum} and results of Section~\ref{sectTHEO} indicates that is not true for higher eigenvalues.\\

\noindent {\bf Theorem \ref{thm:max1:k-com}.}
{\it  The uniform $n$--partition of length $\ell$ is a maximizer for the first
  mixed eigenvalue on the disk among all parts of boundary of length $\ell$ with
  at most $n$ connected components.}

\begin{proof}[Proof of Theorem \ref{thm:maxk}] Following~\cite{denzler:bounds}, we construct a sequence $(\partial_D \Omega)_n$ of Dirichlet's parts such that the corresponding $k$--th eigenfunctions $\phi_{k,n}$ strongly converge to some limit $\phi_k$. By strong convergence we know that the orthogonality relations are preserved and that there's strong convergence $|\phi_{k,n}| \rightarrow |\phi_k|$. Hence, the conclusion follows from the argument used by Denzler~\cite{denzler:bounds}.
\end{proof}

The main idea of the proof of Theorem \ref{thm:max1:k-com} is to use a new kind of rearrangement using the symmetries of the uniform partition. We present here the proof for $k=2$, in order to shorten notation, the general case following from repeating the same process in an appropriate order. We use polar coordinates $(r,\theta)$, on the disk (ie with $\theta \in \bR$ (mod $2\pi$)). The uniform $2$--partition of length $\ell$, denoted $\Gamma_2$, is then represented in Figure~\ref{configQUELC} (ie $\Gamma_2 = \Gamma(\frac{2\pi-\ell}{4})$, see Definition~\ref{betafamily}).

\begin{proof}[Proof of Theorem \ref{thm:max1:k-com} (for $k=2$)] Let $u(r,\theta)$ be a first eigenfunction associated to $\Gamma_2$. We know that this function is symmetric with respect to each isometry of $\Gamma_2$. From this function, for a given $\partial_D \ball \in \mF_{\ell}$ with, at most, two connected components, we build a test function for the associated mixed boundary value problem.
\begin{figure}[h]
\psfrag{p-b}[][][1]{$b$}
\psfrag{a}[][][1]{$a$}
\psfrag{p-a}[][][1]{$\frac{\ell}{2}- a$}
\psfrag{Gab}[][][1]{$\partial_D \ball$}
\psfrag{Gabc}[][][1]{$\partial_N \ball$}
\psfrag{G}[][][1]{$\Gamma_2$}
\psfrag{Gppc}[][][1]{$\Gamma_{2, N}$}
\centering
 \includegraphics{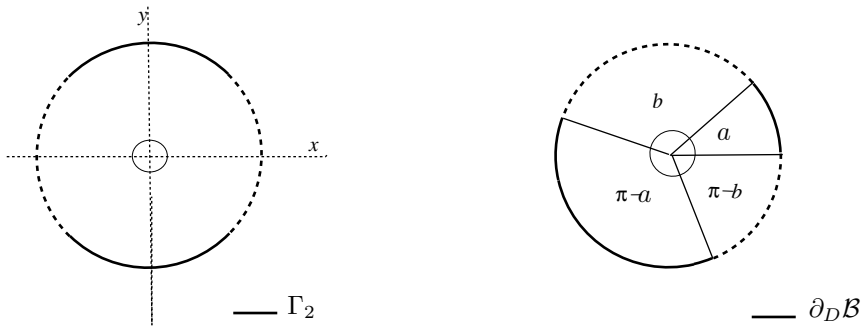}
\caption{Uniform $2$--partition (left) and partition of at most two connected components (right) with $\ell = \pi$.}
\label{configQUELC}
\end{figure}

We denote by $a$ the length of the smallest connected part of $\partial_D \ball$ and by $b$ the length of the smallest connected part of $\partial_N \ball$. Then $a \in [0,\frac{\ell}{2}]$ and $b \in [0,\frac{2\pi-\ell}{2}]$, see Figure~\ref{configQUELC}. We define $$\alpha :=  \frac{1}{2}\left( \frac{2\pi-\ell}{2} - a \right), \;\;\;\; \beta := \frac{1}{2}\left(\frac{2\pi-\ell}{2} - b \right)$$ and

\begin{equation} \widetilde{u}(\theta, r) := \left\{
\begin{array}{lcl}
u(\theta - \alpha + \beta, r) &\;\;& \theta \in [ \alpha , \frac{\pi}{2} - \beta ], \\
u(\theta + \alpha + \beta, r) && \theta \in [\frac{\pi}{2} - \beta , \pi - \alpha], \\
u(\theta + \alpha - \beta, r) && \theta \in [ \pi - \alpha, \frac{3\pi}{2} + \beta], \\
 u(\theta - \alpha - \beta, r) && \theta \in [\frac{3\pi}{2} + \beta, \alpha].
\end{array}
\right.
\end{equation}
It follows from the properties of symmetry and regularity of $u$ that $\widetilde{u}$ is an element of $ H^{1,2}(\ball)$ which is continuous. Moreover the construction ensures that $\widetilde{u}|_{\partial_D \ball} = 0$ and symmetries imply that $\int_{\ball} |\nabla \widetilde{u}|^2 = \int_{\ball} |\nabla u|^2$ and $\int_{\ball} \widetilde{u}^2 = \int_{\ball} u^2$. Hence \begin{equation}
\lambda_1(\partial_D \ball) \leq \frac{\int_{\ball} |\nabla \widetilde{u}|^2}{\int_{\ball} \widetilde{u}^2} = \frac{\int_{\ball} |\nabla u|^2}{\int_{\ball} u^2} = \lambda_1(\Gamma_2).
\end{equation}\end{proof}

\bibliographystyle{abbrv}

\end{document}